\documentclass[11pt]{amsart}

\usepackage{mathpazo}
\usepackage{amsmath,amssymb,amscd,color}
\usepackage{fullpage}
\usepackage{float}
\usepackage{tikz}
\usetikzlibrary{positioning}
\usetikzlibrary{decorations.markings}
\usepackage{comment}

\newtheorem{theorem}{Theorem}[section]
\newtheorem*{theorem*}{Theorem}

\newtheorem{lemma}[theorem]{Lemma}

\newcommand{\teichmuller}{Teichm{\"u}ller{ }}

\definecolor{darkred}{rgb}{0.82,0.0,0.14}

\makeatletter
 \let\c@theorem=\c@subsection
 \let\c@conjecture=\c@subsection
 \let\c@lemma=\c@subsection
 \let\c@proposition=\c@subsection
 \let\c@claim=\c@subsection
 \let\c@question=\c@subsection
 \let\c@criterion=\c@subsection
 \let\c@vfconj=\c@subsection
 \let\c@definition=\c@subsection
 \let\c@notation=\c@subsection
 \let\c@remark=\c@subsection
 \let\c@example=\c@subsection
 \let\c@equation=\c@subsection
 \let\c@figure=\c@subsection
 \let\c@wrapfigure=\c@subsection

\makeatother

\begin{document}

\title{Minimally intersecting filling pairs on the punctured surface of genus two}

\author[Jeffreys]{Luke Jeffreys}
\address{\hskip-\parindent
     School of Mathematics \& Statistics\\
     University of Glasgow\\
     University Place\\
      Glasgow G12 8SQ UK}
\email{l.jeffreys.1@research.gla.ac.uk}


\begin{abstract}
In this short note, we construct a minimally intersecting pair of simple closed curves that fill a genus 2 surface with an odd, greater than 3, number of punctures. This finishes the determination of minimally intersecting filling pairs for all surfaces completing the work of Aougab-Huang and Aougab-Taylor.
\end{abstract}


\maketitle


\section{Introduction}

Let $S_{g,p}$ be an orientable surface of genus $g$ with $p$ punctures. A pair of essential simple closed curves in minimal position on the surface $S_{g,p}$ is said to be a filling pair if the complement of their union is a disjoint collection of disks and once punctured disks. Filling curves on surfaces are key to many fundamental results, such as the Alexander method used to determine if a homeomorphism of a surface is isotopic to the identity, and Thurston's construction of pseudo-Anosov diffeomorphisms.

We define $i_{g,p}$ to be the minimal geometric intersection number of a filling pair on the surface $S_{g,p}$. For closed surfaces, Aougab-Huang calculated the values of $i_{g,0}$ in order to count the number of mapping class group orbits of minimally intersecting filling pairs \cite[Theorem 1.1]{AH}. This count then allowed them to estimate the growth rate of the number of global minima of a topological Morse function they define on the moduli space of Riemann surfaces of genus $g$ \cite[Theorem 1.3]{AH}. 

Aougab-Taylor extended the calculations of $i_{g,p}$ to some cases of $p\neq 0$ \cite[Lemma 3.1]{AT}. This allowed them to construct geodesic rays in the curve graph realising an optimal intersection property \cite[Theorem 1.2]{AT}, answering a question of Margalit. Moreover, filling pairs realising these values of $i_{g,p}$ enabled Aougab-Taylor to construct a certain class of pseudo-Anosov diffeomorphisms that optimise the ratio of \teichmuller to curve graph translation length \cite{AT2}.

The values of $i_{g,p}$ determined so far can be summarised as follows.

\begin{theorem*}[{\cite[Section 2]{AH}}, {\cite[Lemma 3.1]{AT}}]{\label{intro}} The values of $i_{g,p}$ are the following:
\begin{itemize}
\item[(1)] If $g\neq 2,0$ and $p=0$, then $i_{g,p}=2g-1$; 
\item[(2)] If $g\neq 2,0$ and $p\geq1$, then $i_{g,p}=2g+p-2$; 
\item[(3)] If $g=0$ and $p\geq 6$ is even, then $i_{g,p}=p-2$, and if $p\geq{5}$ is odd then $i_{g,p}=p-1$; 
\item[(4)] If $g=2$ and $p\leq 2$, then $i_{g,p}=4$; 
\item[(5)] If $g=2$ and $p\geq 2$ is even, then $i_{g,p}=2g+p-2$. Otherwise, if $p\geq 3$ is odd, then \hbox{$2g+p-2\leq i_{g,p}\leq 2g+p-1$}.
\end{itemize}
\end{theorem*}

In this note, we complete the list by demonstrating the following.

\begin{theorem}\label{result}
Let $g=2$ and $p\geq 3$ be odd, then $i_{g,p}=2g+p-2$.
\end{theorem}

That is, we construct filling pairs on $S_{2,p}$ that realise the lower bound in part (5) of the theorem above. To prove the existence of such filling pairs, we generalise the construction of filling permutations given by Nieland \cite[Theorem 2.1]{MN}, which are themselves generalisations of the filling permutations introduced by Aougab-Huang {\cite[Lemma 2.2]{AH}}. We use these to produce a minimally intersecting filling pair on $S_{2,3}$, and then apply the double-bigon inductive method used by Aougab-Taylor {\cite[Proof of Lemma 3.1]{AT}} to extend to all odd $p\geq 3$.


\section{Preliminaries}

In this section, we will give a short discussion of filling pairs on surfaces and give the generalisation of filling permutations mentioned above. For more details relating to the content of Section \ref{curves}, we refer the reader to the textbook of Farb and Margalit \cite{FM}.

\subsection{Curves, intersection numbers, and filling pairs}\label{curves}

Let $S_{g,p}$ denote the surface of genus $g$ with $p$ punctures. A {\em simple closed curve} $\alpha\subset{S_{g,p}}$ is the image of an embedding of the circle into the surface. We shall say that a simple closed curve $\alpha$ is {\em essential} if it is not homotopic to a point or to a puncture. Equivalently, a simple closed curve is essential if none of the components of its complement are a disk or a once-punctured disk.

Let $\alpha$ be an essential simple closed curve. We denote by $[\alpha]$ the homotopy class of $\alpha$. If $\beta$ is a second essential simple closed curve, then the {\em geometric intersection number} of $\alpha$ and $\beta$, $i(\alpha,\beta)$, is defined to be
\[i(\alpha,\beta):=\min_{\gamma\in[\alpha]}|\gamma\cap\beta|.\]
If the surface is endowed with a complete hyperbolic metric, and $\alpha$ and $\beta$ are chosen to be the unique geodesic representatives in their respective homotopy classes, then we have that $|\alpha\cap\beta|=i(\alpha,\beta)$. We note that the existence of a unique geodesic representative is a consequence of the complete hyperbolic metric. We will say that two essential simple closed curves $\alpha$ and $\beta$ are in {\em minimal position} if $|\alpha\cap\beta|=i(\alpha,\beta)$.

Let $\alpha$ and $\beta$ be a pair of essential simple closed curves in minimal position. We say that $\alpha$ and $\beta$ are a {\em filling pair} if the complement of their union is a disjoint collection of disks and once-punctured disks. Equivalently, they are a filling pair if any other essential simple closed curve intersects at least one of the pair. For a fixed surface $S_{g,p}$, we will let $i_{g,p}$ denote the minimal geometric intersection number of a filling pair on that surface. As discussed above, the values of $i_{g,p}$ were determined in almost all cases in the works of Aougab-Huang \cite{AH} and Aougab-Taylor \cite{AT}. However, in the case of $S_{2,p}$, $p\geq 3$ odd, they showed only the bounds:
\[2g+p-2\leq i_{g,p}\leq 2g+p-1.\]
The lower bound comes from a simple Euler-characteristic argument, realising that the filling pairs form a 4-valent graph on the surface, whereas the upper bound was given by applying a double-bigon construction, detailed in Section \ref{const} below, to a minimally intersecting filling pair on $S_{2,0}$. Our aim is to explicitly construct minimally intersecting filling pairs showing that $i_{g,p}=2g+p-2$.

\subsection{Generalised filling permutations}

We will make use of a generalisation, extending the construction given by Nieland \cite[Theorem 2.1]{MN}, of the filling permutations used by Aougab-Huang in their determinations of $i_{g,0}$ \cite[Lemma 2.2]{AH}. The notation and result that follows is a simple extension of the works of Aougab-Huang and Nieland, however our generalised construction allows us to work with filling pairs on punctured surfaces.

Let $\alpha$ and $\beta$ be a filling pair on the surface $S_{g,p}$ and let $n=i(\alpha,\beta)$. Fix orientations for the curves $\alpha$ and $\beta$ and choose one of the intersection points $x\in\alpha\cap\beta$. Starting at $x$, and following the orientation of $\alpha$, number the arcs of $\alpha$ between consecutive intersection points in order to obtain the set $\{\alpha_{1},\ldots,\alpha_{n}\}$. Similarly, and possibly choosing a different intersection point $y\in\alpha\cap\beta$, construct the set $\{\beta_{1},\ldots,\beta_{n}\}$. Let the set $A=A_{\alpha,\beta}$ be defined by
\[A:=\{\alpha_{1},\beta_{1},\ldots,\alpha_{n},\beta_{n},\alpha_{1}^{-1},\beta_{1}^{-1},\ldots,\alpha_{n}^{-1},\beta_{n}^{-1}\},\]
and identify this set with the set $\{1,2,\ldots,4n\}$.

We now define a {\em filling permutation}, $\sigma=\sigma_{\alpha,\beta}\in\Sigma_{4n}$, associated to the filling pair $\alpha,\beta$, as follows. Firstly, cut $S_{g,p}$ along $\alpha\cup\beta$ to form a collection of $n+2-2g$ polygons with sides labelled by $\alpha$ and $\beta$. Orienting the polygons clockwise, we obtain a labelling of the sides of the polygons by the elements of $A$. We now define $\sigma$ in the following way. If, going around the sides of the polygons in a clockwise direction, the edge labelled by the $j^{th}$ element of $A$ is followed by the edge labelled by the $k^{th}$ element of $A$, then we define $\sigma(j)=k$. We see then that $\sigma$ is an element of the symmetric group $\Sigma_{4n}$ consisting of $n+2-2g$ cycles.

We will also be interested in two more elements of $\Sigma_{4n}$ that have geometric significance. Firstly, we define the permutation $Q=Q_{\alpha,\beta}$ by
\[Q = (1,2,\ldots,4n)^{2n}.\]
Observe that this permutation sends $j$ to $k$ if and only if the $j^{th}$ and $k^{th}$ elements of $A$ are the inverses of one another. Secondly, we define the permutation $\tau=\tau_{\alpha,\beta}$ by
\[\tau = (1,3,5,\ldots,2n-1)(2,4,6,\ldots,2n)(4n-1,4n-3,\ldots,2n+1)(4n,4n-2,\ldots,2n+2).\]
In this case, $\tau$ corresponds to sending an arc of one of the curves to the following arc in the same curve with the same orientation.

Note that we will say that a permutation is parity reversing if it sends odd numbers to even numbers and even numbers to odd numbers.

The following lemma generalises the results of Aougab-Huang \cite[Lemma 2.2]{AH} and Nieland \cite[Theorem 2.1]{MN}, which dealt with the cases $p=0$, $i(\alpha,\beta)=2g-1$, and $p=0$, $i(\alpha,\beta)=n\geq i_{g,0}$, respectively. This amounts to ensuring that any bigons, equivalently 2-cycles of the filling permutation, are dealt with appropriately.

\begin{lemma}\label{perm}
Let $\alpha$ and $\beta$ be a filling pair on $S_{g,p}$ with $i(\alpha,\beta)=n\geq i_{g,p}$. Then $\sigma=\sigma_{\alpha,\beta}$ satisfies the equation
\[\sigma Q\sigma=\tau.\]
Conversely, a parity reversing permutation $\sigma\in\Sigma_{4n}$ consisting of $n+2-2g$ cycles and no more than $p$ number of 2-cycles, and satisfying the above equation, defines a filling pair on $S_{g,p}$ with intersection number $n$.
\end{lemma}

\begin{figure}[H]
\begin{center}
\begin{tikzpicture}[scale = 1.2]
\draw [darkred, line width = 0.3mm, decoration={markings, mark=at position 0.58 with {\arrow{latex}}}, postaction={decorate}] (0,2)--(0,0);
\draw [darkred, line width = 0.3mm, decoration={markings, mark=at position 0.58 with {\arrow{latex}}}, postaction={decorate}] (0,0)--node[left, color = black]{\small{$\sigma(j)^{th}$}} node[right, color = black]{\small{$Q(\sigma(j))^{th}$}} (0,-2);
\draw [blue, decoration={markings, mark=at position 0.58 with {\arrow{latex}}}, postaction={decorate}] (-2,0)--node[below, color = black]{\small{$j^{th}$}}(0,0);
\draw [blue, decoration={markings, mark=at position 0.58 with {\arrow{latex}}}, postaction={decorate}] (0,0)--node[below, color = black]{\small{$\sigma(Q(\sigma(j)))^{th}$}}(2,0);
\end{tikzpicture}
\end{center}
\caption{The filling permutation equation around a vertex.}
\label{vertex}
\end{figure}
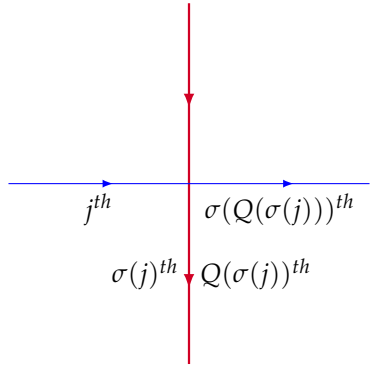

\begin{proof}
Let $j\in\{1,2,\ldots,4n\}$, then the edge labelled by the $j^{th}$ element of $A$ is followed by the edge labelled by the $\sigma(j)^{th}$ element of $A$. As discussed above, we then have that $Q(\sigma(j))$ is the inverse element in $A$ of $\sigma(j)$. Finally, the edge labelled by the $\sigma(Q(\sigma(j)))^{th}$ element of $A$ is the edge following the edge labelled by the $Q(\sigma(j))^{th}$ element of $A$. As can be seen in Figure \ref{vertex}, this is in fact the edge labelled by the arc that follows the arc labelled by $j$ in the same curve. That is, the composition has the same action as the action of $\tau$.

Conversely, suppose that $\sigma\in\Sigma_{4n}$ satisfies the conditions of the lemma. Since $\sigma$ is parity reversing, each of the cycles in $\sigma$ is of even length. Associate to each cycle of $\sigma$ a polygon with the same number of sides. Puncture every 2-gon and then puncture any of the remaining polygons at most once until all $p$ punctures have been placed. This is possible since $n\geq i_{g,p}$ guarantees that $n+2-2g\geq p$. Furthermore, since the number of 2-cycles in $\sigma$ was at most $p$, we do not have any unpunctured bigons that could reduce the intersection number of the resulting curves.

Label each polygon cyclically in a clockwise direction with the elements of the associated cycle of $\sigma$. Now relabel each side with the corresponding element of $A$ and glue the polygons together by gluing each edge to the edge labelled with its inverse. Since every edge occurs once with each orientation, the resulting surface is closed with $p$ punctures.

In the construction so far, we have $n+2-2g$ faces and $2n$ edges, so we have $V-E+F=2-2g$ if and only if we have $n$ equivalence classes of vertices under the gluing of the polygons. Applying the filling permutation equation in Figure \ref{vertex}, we see that four edges will glue together to give a single vertex. In other words, under the gluing, the $4n$ vertices of the polygons form equivalences classes of size 4 and so we have $n$ equivalence classes under the gluing. Hence the resulting surface is $S_{g,p}$.

Finally, another application of the filling permutation equation proves that the $\alpha$-arcs and $\beta$-arcs glue to form a filling pair on $S_{g,p}$ with geometric intersection number equal to $n$.
\end{proof}


\section{Filling pair construction}\label{const}

The filling permutations introduced in the previous section now give us a way to construct filling pairs on punctured surfaces satisfying specific conditions on their intersection number. We now make use of this tool to construct the filling pairs we require to prove Theorem \ref{result}.

\begin{figure}[H]
\begin{center}
\begin{tikzpicture}[scale=2]

\draw (0,0) node{.};
\draw (-2.26,0) node{.};
\draw (2.26,0) node{.};

\draw [darkred, decoration={markings, mark=at position 0.58 with {\arrow{latex}}}, postaction={decorate}] (-2.52,0.26)-- node[above]{$\alpha_{1}$}(-2,0.26);
\draw [blue, decoration={markings, mark=at position 0.58 with {\arrow{latex}}}, postaction={decorate}] (-2,0.26)--node[right]{$\beta_{1}$}(-2,-0.26);
\draw [darkred, decoration={markings, mark=at position 0.58 with {\arrow{latex}}}, postaction={decorate}](-2.52,-0.26)--node[below]{$\alpha_{5}$}(-2,-0.26);
\draw [blue, decoration={markings, mark=at position 0.58 with {\arrow{latex}}}, postaction={decorate}](-2.52,0.26)--node[left]{$\beta_{2}$}(-2.52,-0.26);

\draw [blue, decoration={markings, mark=at position 0.58 with {\arrow{latex}}}, postaction={decorate}](0,1)-- node[above]{$\beta_{4}$}(0.5,0.87);
\draw [darkred, decoration={markings, mark=at position 0.58 with {\arrow{latex}}}, postaction={decorate}](0.87,0.5)--node[above right]{$\alpha_{3}$}(0.5,0.87);
\draw [blue, decoration={markings, mark=at position 0.58 with {\arrow{latex}}}, postaction={decorate}](1,0)--node[right]{$\beta_{3}$}(0.87,0.5);
\draw [darkred, decoration={markings, mark=at position 0.58 with {\arrow{latex}}}, postaction={decorate}](1,0)--node[right]{$\alpha_{5}$}(0.87,-0.5);
\draw [blue, decoration={markings, mark=at position 0.58 with {\arrow{latex}}}, postaction={decorate}](0.87,-0.5)--node[below right]{$\beta_{2}$}(0.5,-0.87);
\draw [darkred, decoration={markings, mark=at position 0.58 with {\arrow{latex}}}, postaction={decorate}](0,-1)--node[below]{$\alpha_{4}$}(0.5,-0.87);
\draw [blue, decoration={markings, mark=at position 0.58 with {\arrow{latex}}}, postaction={decorate}](-0.5,-0.87)--node[below]{$\beta_{4}$}(0,-1);
\draw [darkred, decoration={markings, mark=at position 0.58 with {\arrow{latex}}}, postaction={decorate}](-0.5,-0.87)--node[below left]{$\alpha_{3}$}(-0.87,-0.5);
\draw [blue, decoration={markings, mark=at position 0.58 with {\arrow{latex}}}, postaction={decorate}](-0.87,-0.5)--node[left]{$\beta_{5}$}(-1,0);
\draw [darkred, decoration={markings, mark=at position 0.58 with {\arrow{latex}}}, postaction={decorate}](-0.87,0.5)--node[left]{$\alpha_{1}$}(-1,0);
\draw [blue, decoration={markings, mark=at position 0.58 with {\arrow{latex}}}, postaction={decorate}](-0.5,0.87)--node[above left]{$\beta_{1}$}(-0.87,0.5);
\draw [darkred, decoration={markings, mark=at position 0.58 with {\arrow{latex}}}, postaction={decorate}](-0.5,0.87)--node[above]{$\alpha_{2}$}(0,1);

\draw [blue, decoration={markings, mark=at position 0.58 with {\arrow{latex}}}, postaction={decorate}] (2,0.26)-- node[above]{$\beta_{3}$}(2.52,0.26);
\draw [darkred, decoration={markings, mark=at position 0.58 with {\arrow{latex}}}, postaction={decorate}] (2,-0.26)--node[left]{$\alpha_{4}$}(2,0.26);
\draw [blue, decoration={markings, mark=at position 0.58 with {\arrow{latex}}}, postaction={decorate}](2,-0.26)--node[below]{$\beta_{5}$}(2.52,-0.26);
\draw [darkred, decoration={markings, mark=at position 0.58 with {\arrow{latex}}}, postaction={decorate}](2.52,-0.26)--node[right]{$\alpha_{2}$}(2.52,0.26);

\end{tikzpicture}
\end{center}
\caption{Polygonal decomposition of $S_{2,3}$ associated to the filling permutation $\sigma$.}
\label{poly}
\end{figure}
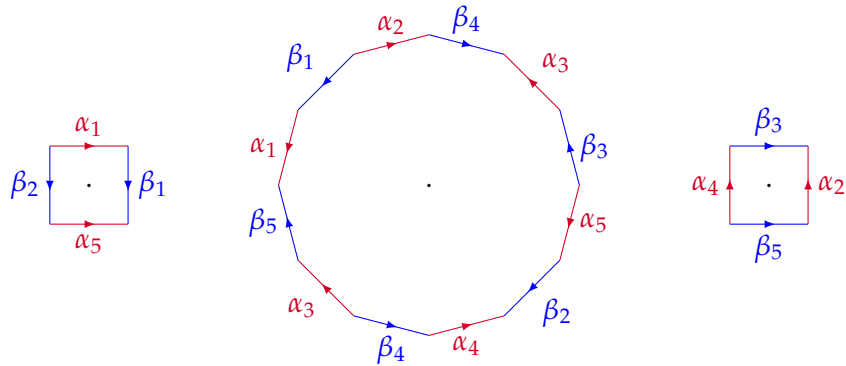

We begin by constructing an appropriate filling pair on the surface $S_{2,3}$. We require that the intersection number of the filling pair is equal to $2g+p-2=5$ and so we must find a permutation $\sigma\in\Sigma_{20}$ consisting of $5+2-2g=3$ cycles and satisfying the conditions of Lemma \ref{perm}. Indeed, the permutation 
\[\sigma=(1,2,19,14)(3,8,15,16,9,4,17,18,5,10,11,12)(6,13,20,7)\]
satisfies the lemma, and the polygonal decomposition of $S_{2,3}$ determined by the associated filling pair is shown in Figure \ref{poly}.

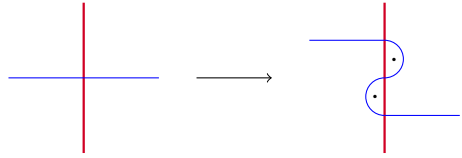
\begin{figure}[H]
\begin{center}
\begin{tikzpicture}[scale = 0.5]
\draw [darkred, line width = 0.3mm] (0,2)--(0,-2);
\draw [blue] (-2,0)--(2,0);
\draw [->] (3,0)--(5,0);
\draw [darkred, line width = 0.3mm] (8,2)--(8,-2);
\draw [blue] (6,1)--(8,1);
\draw [blue] (8,-1) arc (270:90:0.5);
\draw [blue] (8,1) arc (270:90:-0.5);
\draw [blue] (8,-1)--(10,-1);
\draw (8.25,0.5) node{.};
\draw (7.75,-0.5) node{.};
\end{tikzpicture}
\end{center}
\caption{Double-bigon inductive method.}
\label{bigon}
\end{figure}

We now apply the double-bigon inductive method used by Aougab-Taylor {\cite[Proof of Lemma 3.1]{AT}} to complete the proof of Theorem \ref{result}. Namely, suppose that we have a filling pair on the surface $S_{2,p}$ with intersection number equal to $2g+p-2$. Choose an intersection point of the two curves, form two bigons as in Figure \ref{bigon}, and puncture each of these bigons. We now have a filling pair on $S_{2,p+2}$ with intersection number equal to $2g+(p+2)-2$, as required. Since we have constructed a suitable filling pair on $S_{2,3}$, by induction, we have completed the proof of the theorem.

\begin{figure}[H]
\begin{center}
\begin{tikzpicture}
\draw (0,1) .. controls (1,0.9) .. (2,1);
\draw (0,1) .. controls (0.1,0) .. (0,-1);
\draw (0,-1) .. controls (1,-0.9) .. (2,-1);
\draw (2,-1) .. controls (1.9,0) .. (2,1);
\draw [color=darkred] (1,0.95) .. controls (0.9,0) .. (1,-0.95);
\draw [color=darkred, dashed] (1,0.95) .. controls (1.1,0) .. (1,-0.95);
\draw [color=blue] (0.05,0) .. controls (0.8,-0.1) .. (1.95,0);
\draw [color=blue, dashed] (0.05,0) .. controls (1,0.1) .. (1.95,0);
\draw (0,1)node{$\cdot$};
\draw (2,1)node{$\cdot$};
\draw (0,-1)node{$\cdot$};
\draw (2,-1)node{$\cdot$};
\end{tikzpicture}
\end{center}
\caption{A minimally intersecting filling pair on $S_{0,4}$.}
\label{S04}
\end{figure}
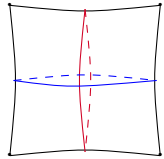

We conclude with the observation that the filling pair in Figure \ref{S04} is a minimally intersecting filling pair on $S_{0,4}$ having intersection number equal to $2=p-2$. As such, the complete list of values of $i_{g,p}$ can be summarised as follows.

\begin{theorem}
The values of $i_{g,p}$ are the following:
\begin{itemize}
\item[(1)] If $g\neq 2,0$ and $p=0$, then $i_{g,p}=2g-1$; 
\item[(2)] If $g\neq 2,0$ and $p\geq1$, then $i_{g,p}=2g+p-2$; 
\item[(3)] If $g=0$ and $p\geq4$, then $i_{g,p}=p-2$ if $p$ is even, and $i_{g,p}=p-1$ if $p$ is odd; 
\item[(4)] If $g=2$ and $p\leq 2$, then $i_{g,p}=4$; 
\item[(5)] If $g=2$ and $p\geq 2$, then $i_{g,p}=2g+p-2$.
\end{itemize}
\end{theorem}

\end{document}